\documentclass{elsart}

\usepackage{amssymb}
\usepackage[mathscr]{eucal}
\usepackage[all]{xy}

\newtheorem{theorem}{Theorem}
\newtheorem{lemma}[equation]{Lemma}
\newtheorem{proposition}[equation]{Proposition}

\newtheorem{conjecture}[equation]{Conjecture}

\renewcommand{\theequation}{\thesection.\arabic{equation}}

\newcommand{\resetL}{\setcounter{equation}{0}} 

\newcommand{\PP}{\mathbb{P}}

\newcommand{\ZZ}{\mathbb{Z}}

\newcommand{\FF}{\mathbb{F}}
\newcommand{\Fq}{\mathbb{F}_q}

\newcommand{\NN}{\mathbb{N}}

\newcommand{\fa}{\mathfrak{a}}
\newcommand{\fA}{\mathfrak{A}}
\newcommand{\fb}{\mathfrak{b}}
\newcommand{\fc}{\mathfrak{c}}
\newcommand{\fd}{\mathfrak{d}}
\newcommand{\fD}{\mathfrak{D}}

\newcommand{\fl}{\mathfrak{l}}

\newcommand{\fm}{\mathfrak{m}}
\newcommand{\fn}{\mathfrak{n}} 
\newcommand{\fp}{\mathfrak{p}}
\newcommand{\fP}{\mathfrak{P}}
\newcommand{\fq}{\mathfrak{q}}

\newcommand{\OO}{\mathcal{O}}

\newcommand{\cN}{\mathcal{N}}
\newcommand{\cM}{\mathcal{M}}

\newcommand{\TT}{\mathcal{T}}

\newcommand{\Tr}{\mathrm{Tr}}

\newcommand{\End}{\mathrm{End}}

\newcommand{\Gal}{\mathrm{Gal}}
\newcommand{\Pic}{\mathrm{Pic}}

\newcommand{\GL}{\mathrm{GL}_2}

\newcommand{\Stab}{\mathrm{Stab}}

\newcommand{\sep}{\mathrm{sep}}

\newcommand{\kinf}{k_\infty}

\newcommand{\Cinf}{{\mathbb{C}}_\infty}

\newcommand{\PGLki}{\mathrm{PGL}_2(k_\infty)}

\newcommand{\uH}{{\underline{H\!\!}\:}}

\newcommand{\Endproof}{\vspace{-11pt}\begin{flushright} $\square$ \end{flushright}}

\renewcommand{\matrix}[4]{\left(\raisebox{-5pt}{$\stackrel{\textstyle
#1 \;\; #2}{#3 \;\; #4}$}\right)}

\journal{Journal of Number Theory}

\begin{document}

\begin{frontmatter}

\title{Higher Heegner points on elliptic curves over 
function fields}
\author{Florian Breuer}
\ead{breuer@math.cts.nthu.edu.tw}
\address{Mathematics Division, National Center for Theoretical Sciences \\
National Tsing-Hua University, Hsinchu, Taiwan}

\begin{abstract}
%
Let $E$ be a modular elliptic curve defined over a rational
function field $k$ of odd characteristic. We construct a sequence
of Heegner points on $E$, defined over a $\ZZ_p^{\infty}$-tower of
finite extensions of $k$, and show that these Heegner points
generate a group of infinite rank.  This is a function field
analogue of a result of C.~Cornut and V.~Vatsal.
\end{abstract}


\end{frontmatter}

\section{Introduction}
\resetL

Heegner points are a way of constructing explicit points of
infinite order on modular elliptic curves, and have been used to
great effect over number fields. More recently, people have
produced analogous constructions on elliptic curves over global
function fields. The first such construction is due to M.L.~Brown
\cite{Brown94}, but there seem to be a number of errors in his
paper. H.-G.~R\"{u}ck and U.~Tipp \cite{RueckTipp} have proved an
analogue of the celebrated Gross-Zagier formula in a special case,
and I.~Longhi \cite{Longhi} and A.~P\'al \cite{Pal} have
(independently) constructed Heegner points of infinite order on
elliptic curves, using the $p$-adic approach due to M.~Bertolini
and H.~Darmon. D.~Ulmer has announced a more general version of
the Gross-Zagier formula, which, combined with his non-vanishing
result for $L$-functions, yields the full Birch and
Swinnerton-Dyer conjecture for elliptic curves over function
fields when their analytic rank is $\leq 1$ and the characteristic
is $p>3$. For a nice survey of all this, see \cite{Ulmer}.

In this paper we prove a function field version of a result of
C.~Cornut \cite{CornutInvent} and V.~Vatsal \cite{Vatsal02} on
higher Heegner points on modular elliptic curves, which had been
conjectured by B.~Mazur \cite{Mazur83}. Their result has some
important consequences, including J.~Nekov\'{a}\v{r}'s celebrated
results concerning the parity of ranks of Selmer groups. See
\cite{CornutInvent} for a discussion.

Cornut actually found two proofs, the second of which
\cite{CornutCRAS} uses a known case of the Andr\'e-Oort conjecture
due to B.~Moonen, and is much simpler. We use the function field
analogue of a special case of the Andr\'e-Oort Conjecture, proved
in \cite{BreuerPrep}, and then follow closely Cornut's second
proof.

Let $k$ be a global function field with field of constants $\Fq$,
where $q$ is a power of the odd prime $p$. Let $E$ be an elliptic
curve defined over $k$ with non-constant $j$-invariant (we say $E$
is non-isotrivial). Then, replacing $k$ by a finite extension if
necessary, we can choose a place $\infty$ of $k$ such that $E$ has
multiplicative reduction at $\infty$. Let $\kinf$ denote the
completion of $k$ at $\infty$, and set $\Cinf =
\hat{\bar{k}}_{\infty}$ the completion of an algebraic closure of
$\kinf$. Furthermore, we let $A$ be the ring of functions in $k$
regular outside $\infty$. It is a Dedekind domain with finite
class number $h=|\Pic(A)| = \deg(\infty)h_k$, where
$h_k=|\Pic^0(k)|$ denotes the class number of $k$. By a Drinfeld
module we will always mean a Drinfeld $A$-module of rank 2,
defined over a subfield of $\Cinf$ (in particular, we deal only
with the case of ``generic'' characteristic).

The conductor of $E$ may be written as $\fn\cdot\infty$, where
$\fn$ is an ideal in $A$. Then, by the work of V.G.~Drinfeld (and
A.~Weil, A.~Grothendieck, H.~Jacquet-R.P.~Langlands, P.~Deligne
and Y.G.~Zarhin) we have a modular parametrization
\begin{equation}\label{param}
\pi : X_0(\fn) \longrightarrow E,
\end{equation}
defined over $k$, where $X_0(\fn)$ is the Drinfeld modular curve
parametrizing isomorphism classes of pairs $(\Phi, \Phi')$ of
Drinfeld modules linked by a cyclic isogeny of degree $\fn$.

Now let $K$ be an imaginary quadratic extension of $k$ (i.e. such
that $\infty$ does not split in $K/k$) with the property that all
primes dividing $\fn$ split in $K$ (this is known as the Heegner
hypothesis). There exist infinitely many such fields. Denote by
$\OO_K$ the integral closure of $A$ in $K$, it contains an ideal
$\cN$ such that $\OO_K/\cN \cong A/\fn$. Let $\fp\subset A$ be a
prime not dividing $\fn$, and let $\OO_n = A + \fp^n\OO_K$ be the
order of conductor $\fp^n$ in $\OO_K$. We set $\cN_n =
\cN\cap\OO_n$, so $\OO_n/\cN_n \cong A/\fn$. Then $\OO_n$ and
$\cN_n^{-1}$, viewed as rank 2 lattices in $\Cinf$, correspond to
a pair of Drinfeld modules $(\Phi^{\OO_n},\Phi^{\cN_n^{-1}})$,
linked by a cyclic isogeny of degree $\fn$. Hence they define a
{\em Heegner point} $x_n \in X_0(\fn)$, which is in fact defined
over the ring class field $K[\fp^n]$ of $\OO_n$, as
$\End(\Phi^{\OO_n}) \cong \End(\Phi^{\cN_n^{-1}}) \cong \OO_n$.

We now set $K[\fp^{\infty}] = \cup_{n\geq 1} K[\fp^n]$. Then $G =
\Gal(K[\fp^{\infty}]/K) \cong \ZZ_p^{\infty}\times G_0$, where
$\ZZ_p^{\infty}$ denotes the product of countably many copies of
$\ZZ_p$ ($p$ being the characteristic of $k$), and $G_0=G_{tors}$
is a finite group (Proposition \ref{tower}). This is in marked
contrast to the number field case, where the analogue of $G$
contains only one copy of $\ZZ_p$, and where $p$ plays the role of
$\fp$, and may be chosen. $G_0$ corresponds to a subfield
$H[\fp^{\infty}]$ satisfying $\Gal(H[\fp^{\infty}]/K) \cong
\ZZ_p^{\infty}$.

Denote by
\[
\Tr_{G_0} : E(K[\fp^{\infty}]) \longrightarrow E(H[\fp^{\infty}]);
\quad x \mapsto \sum_{\sigma\in G_0}x^{\sigma}
\]
the $G_0$-trace on $E$. We define the {\em Heegner point} $y_n =
\Tr_{G_0}(\pi(x_n))\in E(H[\fp^{\infty}])$. Our aim is to prove

\begin{theorem}\label{MainTheorem}
Suppose $k=\Fq(T)$ and $\deg(\infty)=1$. Let $I\subset\NN$ be an
infinite subset. Then the group generated by $\{y_n \;|\; n\in
I\}$ in $E(H[\fp^{\infty}])$ has infinite rank.
\end{theorem}

{\bf Remark.} We have tried to avoid the hypothesis on $k$ and
$\infty$ as far as possible in this paper. It is used twice,
firstly in the proof of Proposition \ref{nonconstant} (but which
should still hold for general $k$), and at the very end of the
proof of Theorem \ref{MainTheorem}, where we invoke an analogue of
the Andr\'e-Oort conjecture which is currently only known in this
case. Once a more general case of this conjecture has been proved
- which is the object of current efforts - Theorem
\ref{MainTheorem} should become true for general $k$ and $\infty$.

The layout of this paper is as follows. In \S\ref{SectTower} we
describe the group $\Gal(K[\fp^{\infty}]/K)$, and show that
$E(K[\fp^{\infty}])$ has finite torsion. In \S\ref{SectModular} we
describe the map (\ref{param}) in more detail, and construct a
family of new modular parametrizations in \S\ref{SectDegen} by
means of degeneracy maps between Drinfeld modular curves. Then in
\S\ref{SectFactor} we describe a canonical factorization of cyclic
isogenies between CM Drinfeld modules, which we will use in
\S\ref{SectGeometric} to characterize the geometric action of
Galois on Heegner points. Finally, we deduce Theorem
\ref{MainTheorem} from the Andr\'e-Oort conjecture in
\S\ref{SectProof}.

{\bf Acknowledgements.} The author would like to thank
Ernst-Ulrich Gekeler for help with Proposition \ref{nonconstant},
and Christophe Cornut for showing him that the $y_n$'s actually
generate a group of infinite rank.

\section{The class field tower}\label{SectTower}
\resetL

\begin{proposition}\label{tower}
$G=\Gal(K[\fp^{\infty}]/K) \cong \ZZ_p^{\infty}\times G_0$, where
$G_0 = G_{tors}$ is finite.
\end{proposition}

{\bf Proof.} First, notice that
\begin{equation}\label{limit}
G = \Gal(K[\fp^{\infty}]/K) =
\lim_{\stackrel{\textstyle\longleftarrow}{n}}\Gal(K[\fp^n]/K)
\cong \lim_{\stackrel{\textstyle\longleftarrow}{n}}\Pic(\OO_n).
\end{equation}
Secondly, we have an exact sequence (see e.g. \cite{Neukirch}, \S
I.12)
\[
1 \rightarrow \OO_K^*/\OO_n^* \rightarrow
(\OO_K/\fp^n\OO_K)^*/(\OO_n/\fp^n\OO_K)^* \rightarrow \Pic(\OO_n)
\rightarrow \Pic(\OO_K) \rightarrow 1.
\]
As $\OO_K^*/\OO_n^*$ and $\Pic(\OO_K)$ are bounded, it remains to
examine the behavior of
$(\OO_K/\fp^n\OO_K)^*/(\OO_n/\fp^n\OO_K)^*$ as
$n\rightarrow\infty$. We insert it into the following diagram,
with exact rows and columns:

\[
\xymatrix@R=15pt@M=5pt{ & 1\ar[d] & 1\ar[d] & 1\ar[d] & \\
1 \ar[r] & \frac{\textstyle 1+\fp\OO_n}{\textstyle
1+\fp^n\OO_n}\ar[r]\ar[d] & (\OO_n/\fp^n\OO_n)^*\ar[r]\ar[d] &
(\OO_n/\fp\OO_n)^*\ar[d]\ar[r]
& 1 \\
1\ar[r] & \frac{\textstyle 1+\fp\OO_K}{\textstyle
1+\fp^n\OO_K}\ar[r]\ar[d] & (\OO_K/\fp^n\OO_K)^*\ar[r]\ar[d] &
(\OO_K/\fp\OO_K)^*\ar[r]\ar[d]
& 1 \\
1\ar[r] & \frac{\textstyle 1+\fp\OO_K}{\textstyle
1+\fp\OO_n}\ar[r]\ar[d] & \frac{\textstyle
(\OO_K/\fp^n\OO_K)^*\!\!\!\!}{\textstyle
(\OO_n/\fp^n\OO_K)^*\!\!\!\!} \;\;\ar[r]\ar[d] & \frac{\textstyle
(\OO_K/\fp\OO_K)^*\!\!\!\!}{\textstyle (\OO_n/\fp\OO_n)^*\!\!\!\!}
\;\;\ar[r]\ar[d] & 1 \\
& 1 & 1 & 1 &}
\]
Notice that $\fp^n\OO_K = \fp^n\OO_n$.
As $(\OO_K/\fp\OO_K)^*/(\OO_n/\fp\OO_n)^*$ 
is bounded, we are lead to studying the (multiplicative) group
$H_n := (1+\fp\OO_K)/(1+\fp\OO_n)$.

By computing the cardinality of various groups in the diagram, one
finds that $H_n$ is a $p$-group, of order $|\fp|^{n-1}$. Let $x
\in 1 + \fp\OO_K$, we will examine its order in $H_n$. Let
$s=\lceil\log_p(n+1)\rceil$, then we find that $x^{p^s} \in
1+\fp^{n+1}\OO_K \subset 1 + \fp\OO_n$. It follows that $H_n$ is
annihilated by $p^s$, hence the number of generators of $H_n$ is
at least $\log_p(|\fp|)(n-1)/\lceil\log_p(n+1)\rceil
\rightarrow\infty$ as $n\rightarrow\infty$. On the other hand,
suppose the order of $x$ in $H_n$ is bounded independently of $n$,
say $x^{p^r} \in 1+\fp\OO_n$ for all $n$. But then $x^{p^r} \in
\cap_{n=1}^{\infty}(1+\fp\OO_n) \subset A$, and so $x\in A$ to
begin with (recall that $K = k(\sqrt{D})$ for some square-free
$D\in A$, and that $p$ is odd), so $x\in A\cap(1+\fp\OO_K) \subset
1 + \fp\OO_n$. We have shown that
${\displaystyle\lim_{\longleftarrow}}\: H_n \cong \ZZ_p^{\infty}$
and the proposition now follows. \Endproof

The following result will be crucial.

\begin{lemma}\label{torsion}
$E_{tors}(K[\fp^{\infty}])$ is finite.
\end{lemma}

{\bf Proof.} Let $\fl\nmid \fp\fn$ be a prime of $k$ which is
inert and principal in $K$. Then $E$ has good reduction at $\fl$,
and $\fl$ splits completely in $K[\fp^n]$, hence the residue field
of $K[\fp^n]$ at $\fl$ is just $\FF_{\fl}=\OO_K/\fl$ for every
$n\geq 0$. It follows that reduction mod $\fl$ induces an
injection of the prime-to-$p$ part of $E_{tors}(K[\fp^{\infty}])$
into $\tilde{E}(\FF_{\fl})$, which is finite.

Let $K^{\sep}$ denote the separable closure of $K$. We complete
the proof by showing that $E(K^{\sep})\cap E[p^{\infty}]$ is
finite, for any non-isotrivial elliptic curve $E/K$. Indeed, if
$E[p^n]\subset E(K)$, then $j=j(E)$ is a $p^n$\/th power in $K$,
as can be seen by factoring the multiplication by $p^n$-map into
$[p^n]=f\circ g$, with $\ker(g)=E[p^n]$ and $f$ the $p^n$\/th
power Frobenius. Now, let $n$ be such that $j$ is not a $p^n$\/th
power in $K$. Then $j$ is a $p^n$\/th power in $K'=K(E[p^n])$, and
it follows that $K'/K$ is not separable. Alternatively, one may
apply a general criterion of J.F.~Voloch for abelian varieties
(\cite{Voloch95}, \S4). \Endproof

\section{Modular parametrizations}\label{SectModular}
\resetL

As the literature already contains excellent expositions of the
theory of Drinfeld modular curves and the parametrization
(\ref{param}), such as \cite{GekelerDMC}, \cite{Gekeler95},
\cite{GekelerReversat}, and \cite{Longhi} (\S1.4 and \S1.5), we
will not attempt a detailed account. Instead, we only recall here
some of the results and notations that we will need, and refer the
reader to \cite{GekelerReversat} for the details.

Every projective $A$-module of rank 2 is isomorphic to $Y_{\fa} =
A\times\fa \subset k^2$ for an ideal $\fa\subset A$, so in
particular their isomorphism classes correspond to $\Pic(A)$. The
group $\GL(k)$ acts on $k^2$ from the right. Let
$x=[\fa]\in\Pic(A)$ and let $\Gamma_x = \Stab_{\GL(k)}(Y_{\fa})$
be the stabilizer of $Y_{\fa}$. Denote by $\Omega =
\Cinf\smallsetminus\kinf$ the Drinfeld upper half-plane, and
choose $z\in\Omega$. We map $k^2$ into $\Omega$ by sending $(a,b)$
to $az+b$. Then the image of any projective rank 2 $A$-module
under this map is a {\em lattice} in $\Cinf$, i.e. a discrete
projective $A$-submodule of $\Cinf$ of rank 2. $\GL(k)$ acts on
$\Omega$ from the left by fractional linear transformations, and
on the lattices this action corresponds to the right action on
$k^2$.

Let $M_0(\fn)$ be the coarse moduli scheme for the moduli problem
``pairs of Drinfeld modules linked by a cyclic isogeny of degree
$\fn$''. This is equivalent to the problem ``pairs of lattices
$\Lambda_1 \subset_{\fn} \Lambda_2$'', where the $\subset_{\fn}$
notation means that $\Lambda_2/\Lambda_1 \cong A/\fn$ as
$A$-modules. We denote by $Y_0(\fn) = M_0(\fn)\times k$ the base
extension to $k$, and by $X_0(\fn)$ the smooth projective model of
$Y_0(\fn)$, which may be obtained from $Y_0(\fn)$ by including
finitely many cusps. The curve $X_0(\fn)$ has $h=|\Pic(A)|$
irreducible components, each defined over the Hilbert class field
$H$ of $(k,A)$, and which we denote by $X_x$ for $x=[\fa]
\in\Pic(A)$. The group $\Gal(H/k)$ permutes the components by
$(\fa,H/k)X_{[\fb]} = X_{[\fa^{-1}\fb]}$.

Each component $X_x$ has an analytic parametrization
\[
X_x(\Cinf) \cong \Gamma_x(\fn)\diagdown\Omega^*, \qquad\mbox{where
$\Omega^* = \Omega\cup\PP^1(k)$, and}
\]
\[
\Gamma_x(\fn) = \left\{\matrix{a}{b}{c}{d} \in \Gamma_x \;|\; c
\equiv 0 \bmod \fn \right\}.
\]

Let $E/k$ be a non-isotrivial elliptic curve with conductor
$\fn\cdot\infty$, where $\fn\subset A$. Then there is a surjective
morphism defined over $k$
\begin{equation}\label{param2}
\pi : X_0(\fn) \longrightarrow E.
\end{equation}
We may suppose that $E$ is maximal in its $k$-isogeny class with
respect to the map (\ref{param2}), we call it the strong Weil
curve. Then this map can be given explicitly on the $\Cinf$-valued
points, which we will describe next.

Let $\TT$ be the Bruhat-Tits tree for $\PGLki$, and denote by
$\uH(\TT,\ZZ)$ the group of harmonic cochains on $\TT$ with values
in $\ZZ$, and define the subgroup $\uH_!(\TT,\ZZ)^{\Gamma_x(\fn)}$
of cochains invariant under $\Gamma_x(\fn)$ and with compact
support on $\Gamma_x(\fn)\diagdown\TT$. Then these cochains
correspond to certain automorphic forms on $\GL$. Moreover, to $E$
one associates, for each $x=[\fa]\in\Pic(A)$, a primitive Hecke
newform $\varphi_x\in\uH_!(\TT,\ZZ)^{\Gamma_x(\fn)}$. To
$\varphi_x$ we associate furthermore a holomorphic theta function
$u_x : \Omega \rightarrow \Cinf^*$ with multiplier $c_x$, (i.e.
$u_x(\alpha z) = c_x(\alpha)u_x(z)$ for all
$\alpha\in\Gamma_x(\fn)$). Let $\Delta_x = \{c_x(\alpha) \;|\;
\alpha\in\Gamma_x(\fn)\}$, which is a multiplicative lattice in
$\Cinf^*$. Then $E$, which has multiplicative reduction at
$\infty$, is isomorphic over $\kinf$ to the Tate curve
$\Cinf^*/\Delta_x$.
We have the explicit parametrization
\[
\xymatrix@R-5pt@M+2pt{X_x(\Cinf)\ar[r]^{\pi_x} & E(\Cinf) \\
\qquad\qquad Y_x(\Cinf) = \Gamma_x(\fn)\diagdown\Omega%
\ar@{^{(}->}[u]\ar[r]^(.7){u_x}
& \Cinf^*/\Delta_x\ar[u]_{\cong} \\
[z]\ar@{|->}[r] & u_x(z) \bmod \Delta_x.}
\]

\section{Degeneracy maps}\label{SectDegen}
\resetL

We next define degeneracy maps between Drinfeld modular curves.
Let $\fm\subset A$ be an ideal coprime to $\fn$. A generic point
of $X_0(\fm\fn)$ can be written as $(\Phi,\Phi/C)$, where $\Phi$
is a Drinfeld module, and $C \cong A/\fm\fn$ an $A$-submodule of
$\Phi_{tors}$. For any divisor $\fd|\fm\fn$ we denote by $C[\fd]
\cong A/\fd$ the $\fd$-torsion submodule of $C$. Now, for every
$\fd | \fm$ we define the $\fd$th degeneracy map
\begin{eqnarray*}
\delta_{\fd} : X_0(\fm\fn) & \longrightarrow & X_0(\fn) \\
(\Phi,\Phi/C) & \longmapsto & (\Phi/C[\fd],\Phi/C[\fd\fn]),
\end{eqnarray*}
which maps the $[\fa]$-component of $X_0(\fm\fn)$ to the
$[\fd^{-1}\fa]$-component of $X_0(\fn)$. In this way, we may
define the Hecke correspondence $T_{\fm} \subset X_0(\fn)^2$ as
the image of $X_0(\fm\fn)$ under the map
$\delta_1\times\delta_{\fm}$.

More generally, let $\tau(\fm)$ be the number of divisors of
$\fm$, then we define the full degeneracy map $\delta :
X_0(\fm\fn) \rightarrow X_0(\fn)^{\tau(\fm)}$ as the product
$\delta = \prod_{\fd|\fm} \delta_{\fd}$. Composing with $\pi$, we
obtain a new parametrization of $E$ by $X_0(\fm\fn)$:
\begin{equation}\label{newparam}
\begin{array}{rccccccc}
\pi': & X_0(\fm\fn) & \stackrel{\delta}{\longrightarrow} &
X_0(\fn)^{\tau(\fm)} & \stackrel{\pi}{\longrightarrow} &
E^{\tau(\fm)} & \stackrel{\Sigma}{\longrightarrow} & E \\
& z & \longmapsto & (z_{\fd})_{\fd|\fm}  & \longmapsto &
(\pi(z_{\fd}))_{\fd|\fm} & \longmapsto &
\sum_{\fd|\fm}\pi(z_{\fd}).
\end{array}
\end{equation}

\begin{proposition}\label{nonconstant}
The morphism $\pi' : X_0(\fm\fn) \rightarrow E$ is defined over
$k$. Suppose that $k=\Fq(T)$ and $\deg(\infty)=1$. Then $\pi' :
X_0(\fm\fn) \rightarrow E$ is surjective.
\end{proposition}

{\bf Proof.} It is clear that $\Gal(H/k)$ leaves $\pi'$ invariant,
and hence $\pi'$ is defined over $k$.

Now we suppose that $k=\Fq(T)$ and $\deg(\infty)=1$. Then, after
replacing $T$ by another generator if necessary, we may take
$A=\Fq[T]$. In particular, $\Pic(A)$ is trivial, and the modular
curves $X_0(\fm\fn)$ and $X_0(\fn)$ are irreducible and defined
over $k$. Analytically, the map $\pi'$ is given by
\[
\pi'([z]) = \prod_{\fd|\fm}u(dz) \bmod \Delta \;\;\; \in
\Cinf^*/\Delta = E(\Cinf),
\]
where $[z]$ denotes the class of $z\in\Omega$ in $Y_0(\fm\fn)$, $u
= u_{\varphi}$ is the theta function associated to the newform
$\varphi$, and $d\in A$ is the monic generator of the ideal $\fd$,
for each $\fd|\fm$. We need to show that the map $u'(z) =
\prod_{\fd|\fm}u(dz)$ is not constant.


Denote by $\OO_{\Omega}(\Omega)$ the ring of rigid holomorphic
functions on $\Omega$, then there is an exact sequence
\[
1 \longrightarrow \Cinf^* \longrightarrow \OO_{\Omega}(\Omega)^*
\stackrel{r}{\longrightarrow} \uH(\TT,\ZZ) \longrightarrow 0.
\]
Furthermore, for any $f\in\GL(k)$ we have $r(u \circ f) = \varphi
\circ f \in \uH(\TT,\ZZ)$, so if $u'$ is constant, then we get
\begin{equation}\label{constanteq}
0 = r(u) = \varphi + \sum_{\fd|\fm,\;\fd \neq 1} \varphi\circ
d\qquad \in\uH(\TT,\ZZ),
\end{equation}
where $d$ denotes the matrix $\matrix{d}{0}{0}{1}$. Now we examine
the Fourier coefficients of these terms. From Proposition 2.10 of
\cite{Gekeler95a} follows that the ``first'' Fourier coefficient
of $\varphi \circ d$ is $c(\varphi\circ d, (1)) = 0$ if $d\not\in
\Fq^*$. On the other hand, by a result of Atkin-Lehner (see \S1.4
of \cite{Gekeler85}) follows that every newform $\varphi\neq 0$
satisfies $c(\varphi,(1)) \neq 0$). Thus (\ref{constanteq})
implies that $\varphi=0$, a contradiction. \Endproof

{\bf Remark.} In the case of general $k$ and $\infty$, Proposition
\ref{nonconstant} should still hold. Indeed, for each component
$X_x$ of $X_0(\fm\fn)$, one has a similar analytical description
of $\pi'$ and an equality of the form (\ref{constanteq}). But
these now involve a combination of different theta functions
$u_y$, one for each component $X_y$ of $X_0(\fn)$ such that
$y=x[\fd^{-1}]$ in $\Pic(A)$ for some $\fd|\fm$. This may be
simplified by considering only those $\fm$ for which every
$\fd|\fm$ is principal (in our application in \S\ref{SectProof}
this amounts to choosing $K$ such that only principal primes of
$k$ ramify in $K$), but one must still calculate the Fourier
coefficients in this case.

\section{Canonical factorizations of isogenies}\label{SectFactor}
\resetL

In this section, we will describe the following canonical
factorization of cyclic isogenies between CM Drinfeld modules.

\begin{proposition}\label{factorlemma2}
Let $f : \Psi_1 \rightarrow \Psi_2$ be a cyclic isogeny of degree
$\fd$ between Drinfeld modules with complex multiplication by
orders $\End(\Phi_i)$ of conductor $\fc_i$ in $K$, for $i=1,2$.
Then we have a commutative diagram of cyclic isogenies
\[
\xymatrix{\Psi_1 \ar[r]^f\ar[d]^{f_1} & \Psi_2\ar[d]^{f_2} \\
\Psi'_1 \ar[r]^{f'} & \Psi'_2} \qquad\qquad\qquad\xymatrix{
\fc_1\ar[r]^{\fd}\ar[d]^{\fd_1} & \fc_2\ar[d]^{\fd_2} \\
\fc\ar[r]^{\fd'} & \fc. }
\]
Here $\End(\Psi'_1)=\End(\Psi'_2)=\OO_{\fc}$ is an order of
conductor $\fc$ in $K$, and $f,f_1,f_2,f'$ are of degree
$\fd,\fd_1,\fd_2,\fd'$, respectively. This data is summarized in
the right-hand diagram. Furthermore, we have $\fd=\fd_1\fd_2\fd'$,
$\fc=\fc_1/\fd_1=\fc_2/\fd_2$ and $\fc + \fd' = A$. Lastly, there
is an ideal $\fD \subset \OO_K$ with $\OO_K/\fD \cong A/\fd'$ such
that $f'$ corresponds to $\fD\cap\OO_{\fc}$.
\end{proposition}

{\bf Proof.} We follow closely the appendix of \cite{CornutCRAS}.
We will use the following notation. Let $\fa \subset \fb$ be
lattices (all lattices here have rank $2$ and are contained in
$K\subset\Cinf$), then we write $\fa \subset_{\fm} \fb$ if
$\fb/\fa\cong A/\fm$. Also, the conductor $c(\fa)$ of the lattice
$\fa$ is the conductor of the order $\End(\fa)$ in $\OO_K$.

Suppose $\fa \subset_{\fd} \fb$. Following \S9 of
\cite{CornutCRAS} word for word, one arrives at an inclusion of
lattices
\begin{equation}
\fa \subset_{\fd_1} \OO_{\fc}\fa \subset_{\fd'} \fd_2\OO_{\fc}\fb
\subset_{\fd_2} \fb.
\end{equation}
Here $\OO_{\fc}$ is an order (of conductor $\fc$), which is
maximal with respect to the property $\OO_{\fc}\fa \subset \fb$.
Let $\fc_1=c(\fa), \fc_2=c(\fb)$, then $\fd_1 = \fc_1/\fc$,
$\fd_2=\fc_2/\fc$ and $\fd'=\fd/\fd_1\fd_2$. To continue the
argument, we need the following folklore result.

\begin{lemma}\label{sublattice}
Let $\fb$ be a lattice of conductor $\fc$, and let $\fq\subset A$
be a prime dividing $\fc$. Then there exactly $|\fq|+1$
sublattices $\fa \subset_{\fq} \fb$. One of them (given by
$\fa=\fq\OO_{\fc/\fq}\fb$) has conductor $c(\fa)=\fc/\fq$ and the
other $|\fq|$ lattices have conductor $c(\fa)=\fc\fq$.
%
\Endproof
\end{lemma}

We want to show that $\fc$ and $\fd'$ are coprime. Let $\fq\subset
A$ be a prime. If $\fq | \fc$, then by Lemma \ref{sublattice},
$\fa'=\OO_{\fc/\fq}\fa$ is the unique $\OO_{\fc}$-stable lattice
satisfying $\fq\fa' \subset_{\fq} \OO_{\fc}\fa$, or, equivalently,
$\OO_{\fc}\fa \subset_{\fq} \fa'$. Now suppose $\fq$ divides
$\fd'$. Then $\fd_2\OO_{\fc}\fb/\OO_{\fc}\fa \cong A/\fd'$ has a
unique $A$-submodule isomorphic to $A/\fq$, namely its
$\fq$-torsion submodule: $(\fd_2\OO_{\fc}\fb/\OO_{\fc}\fa)[\fq] =
(\fd_2\OO_{\fc}\fb\cap\fq^{-1}\OO_{\fc}\fa)/\OO_{\fc}\fa$, hence
$\fa''=\fd_2\OO_{\fc}\fb\cap\fq^{-1}\OO_{\fc}\fa$ is also
$\OO_{\fc}$-stable and satisfies $\OO_{\fc}\fa \subset_{\fq}
\fa''$, and is furthermore contained in $\fb$. Thus, if $\fq$
divides both $\fc$ and $\fd'$, then $\OO_{\fc/\fq}\fa = \fa' =
\fa'' \subset \fb$, which contradicts the maximality of
$\OO_{\fc}$. So $\fc$ and $\fd'$ are coprime.

Lastly, we set $\fD = \fd_2^{-1}\OO_K\fb^{-1}\fa$. Then we see
easily that $\fD \subset_{\fd'} \OO_K$. Also, $\fD_{\fc} =
\fD\cap\OO_{\fc}$ is invertible in $\OO_{\fc}$ and we have
$\fd_2\OO_{\fc}\fb = \fD_{\fc}^{-1}\OO_{\fc}\fa$. Now, using the
equivalence between lattices and Drinfeld modules, Proposition
\ref{factorlemma2} follows. \Endproof

\section{Geometric action}\label{SectGeometric}
\resetL

Let $\Phi_n = \Phi^{\OO_n}$. We say an element $\sigma\in G$ is
``geometric'' if there exists a cyclic isogeny $f_{\fd} : \Phi_n
\rightarrow \Phi_n^{\sigma}$ of fixed degree $\fd$ for infinitely
many $n\in\NN$. In particular, if $\fd + \fn = A$, then each
$(x_n,x_n^{\sigma})$ lies in the Hecke correspondence $T_{\fd}
\subset X_0(\fn)^2$.

Let $\fp_1,\ldots,\fp_g$ be the primes $\neq\fp$ of $k$ which
ramify in $K/k$, and let $\fP_1,\ldots,\fP_g$ be the primes of $K$
lying above them. Denote by $\sigma_i = (\fP_i,K[\fp^{\infty}]/K)$
the Frobenius elements, and by $G_1 = \langle
\sigma_1,\ldots,\sigma_g \rangle$ the group they generate. Suppose
$[\fp_i]$ has order $e_i$ in $\Pic(A)$. Then $\sigma_i$ has order
$2e_i$ in $G$, and it follows that $G_1 \subset G_0$. Let
$\fm=\fp_1^{2e_1-1}\cdots\fp_g^{2e_g-1}$, which is coprime to
$\fn$ (recall that every prime in $\fn$ splits in $K/k$). Then the
elements of $G_1$ are in a one-to-one correspondence with the
divisors $\fd$ of $\fm$, via $\fd=\prod_{i=1}^{g}\fp_i^{n_i}
\mapsto \sigma_{\fd}= \prod_{i=1}^{g}\sigma_i^{n_i}$.

\begin{proposition}\label{geometric}
$G_1$ is the subgroup of geometric elements of $G_0$. More
precisely,
\begin{enumerate}

\item Let $\sigma=\sigma_{\fd}\in G_1$ for some $\fd|\fm$. Then
$(x_n,x_n^{\sigma})\in T_{\fd}\subset X_0(\fn)^2$ for all
$n\in\NN$. In particular, $(x_n^{\sigma})_{\sigma\in G_1} \subset
X_0(\fn)^{|G_1|}$ lies in the image of $X_0(\fm\fn)$ under the
full degeneracy map $\delta$ described in \S\ref{SectDegen}.

\item Conversely, let $\sigma\in G_0$ and suppose that there exist
cyclic isogenies $f_n :\nobreak \Phi_n \rightarrow
\Phi_n^{\sigma}$ of fixed degree $\fD_{\sigma}$ for infinitely
many $n\in\NN$. Then $\sigma\in G_1$.

\end{enumerate}
\end{proposition}

{\bf Proof.} From (\ref{limit}) follows that each $\sigma\in G$
can be written in the form $\sigma = (\sigma_1, \sigma_2,
\ldots)$, where for each $n\in\NN$ $\sigma_n\in\Gal(K[\fp^n]/K)$
corresponds to an invertible ideal $\fA_n = \fA\cap\OO_n$ in
$\OO_n$, for some $\fA \subset\OO_K$.

For each $m\geq n$ the theory of complex multiplication gives an
isogeny $f_m : \Phi_n \rightarrow \Phi_n^{\sigma_m}$ with
$\ker(f_m)\cong\OO_m/\fA_m$. Now let $\fa\subset A$ lie under
$\fA=\fA_n\OO_K$. If $\sigma\in G_1$, then $\fa$ is a product of
primes which ramify in $K/k$. Thus each
$\ker(f_m)\cong\OO_m/\fA_m\cong A/\fa$ is cyclic. It follows that
$(\Phi_n, \Phi_n^{\sigma})$ lies on the curve $Y_0(\fa)$.
Furthermore, $\fa$ is prime to $\fn$, hence $\sigma$ is compatible
with $\fn$-isogenies, and it follows that the pair of Heegner
points $(x_n, x_n^{\sigma})$ lies on the Hecke correspondence
$T_{\fa} \subset X_0(\fn)^2$. This proves part (1) of Proposition
\ref{geometric}.

To prove part (2), we let $I \subset \NN$ be an infinite subset,
$\sigma\in G$, and we suppose that there exists a cyclic isogeny
$f_n : \Phi_n \rightarrow \Phi_n^{\sigma}$ of degree $\fd$ for all
$n\in I$. One potential source of trouble is the fact that $\fp$
might divide $\fd$. Write $\fd=\fp^t \fd'$ with $\fp\nmid \fd'$.
We apply Proposition \ref{factorlemma2} to the case where $\Psi_1
= \Phi_n$ and $\Psi_2 = \Phi_n^{\sigma}$, and $n \in I$ satisfies
$n\geq t/2$. Then we see that $t=2r$ is even, and we have a
commutative diagram (with the relevant conductors and degrees
shown on the right)
\[
\xymatrix{\Psi_1 \ar[r]^f\ar[d]^{f_1} &
\Psi_2\ar[d]^{f_2}\ar@{=}[r] & \Psi_1^{\sigma}\ar[d]^{f_1^{\sigma}} \\
\Psi'_1 \ar[r]^{f'} & \Psi'_2\ar@{=}[r]^{?} & {\Psi'_1}^{\sigma}}
\qquad\qquad\qquad\xymatrix{ \fp^n\ar[r]^{\fd}\ar[d]^{\fp^r} &
\fp^n\ar[d]^{\fp^r}\ar@{=}[r] & \fp^n\ar[d]^{\fp^r} \\
\fp^{n-r}\ar[r]^{\fd'} & \fp^{n-r} & \fp^{n-r} }
\]
We claim that $\Psi'_2 = {\Psi'_1}^{\sigma}$. Indeed, we have maps
$f_1^{\sigma} : \Psi_2 = \Psi_1^{\sigma} \rightarrow
{\Psi'_1}^{\sigma}$, and $f_2 : \Psi_2 \rightarrow \Psi'_2$, both
of which correspond to lattice inclusions of the form
$\Lambda_1,\Lambda_2 \subset_{\fp^r} \Lambda$, where $c(\Lambda_1)
= c(\Lambda_2) = \fp^{n-r} = c(\Lambda)/\fp^{r}$. It follows from
Lemma \ref{sublattice} (and induction on $r$) that $\Lambda_1 =
\Lambda_2$, which proves the claim.

Now we restrict our attention to $f' : \Psi'_1 \rightarrow
{\Psi'_1}^{\sigma}$. As $f'$ has degree $\fd'$, which is prime to
$\fp$, it follows that $f'$ corresponds, from Proposition
\ref{factorlemma2}, to an invertible ideal $\fD_n\subset\OO_K$
such that $\OO_K/\fD_n \cong A/\fd'$ for all $n$. This leaves only
finitely many possibilities for $\fD_n$, hence, by restricting
$I$, we may assume $\fD_n = \fD \subset \OO_K$ for all $n\in I$.
It now follows that $(\fD,K[\fp^{\infty}]/K) = \sigma$. In
particular, if $\sigma\in G_0$, then $\sigma$ has finite order,
and there exists an integer $s$ such that $\fD^s$ is principal in
$\OO_K$. Moreover, $\fD^s\cap\OO_n$ is principal for each $n\in
I$, and so $\fD^s$ is generated by an element $d\in\cap_{n\in
I}\OO_n = A$. Denote by $z \mapsto \bar{z}$ the non-trivial
element of $\Gal(K/k)$. Then for every prime $\fP$ of $\OO_K$
dividing $\fD$ we see that $\bar{\fP}$ also divides $\fD$. But as
$\OO_K/\fD \cong A/\fd'$ is cyclic, this is only possible if $\fP$
is ramified in $K/k$, and hence $\sigma\in G_1$. This concludes
the proof of Proposition \ref{geometric}. \Endproof

{\bf Remark.} We have also shown that the geometric elements of
$G$ are the elements of the form $(\fD,K[\fp^{\infty}]/K)$ for
some $\fD\subset\OO_K$ with $\OO_K/\fD$ cyclic. Thus they form a
countable subgroup of $G$.

\section{Proof of Theorem \ref{MainTheorem}}\label{SectProof}
\resetL

Lemma \ref{torsion} says that $E_{tors}(K[\fp^{\infty}])$ is
finite, so in particular, all these torsion points are defined
over a finite extension of $k$. Now Theorem \ref{MainTheorem} will
follow if the fields of definition of the $y_n$'s grow with $n$.
We prove this as follows.

Let $\fm = \fp_1\fp_2\cdots\fp_g \subset A$ as in
\S\ref{SectGeometric}, and recall the full degeneracy map $\delta
: X_0(\fm\fn) \rightarrow X_0(\fn)^{|G_1|}$ from
\S\ref{SectDegen}. From Proposition \ref{geometric} follows that
the point $(x_n^{\sigma})_{\sigma\in G_1}$ lies in the image of
$\delta$, and we denote by $x'_n\in X_0(\fm\fn)(K[\fp^{\infty}])$
its preimage, which is given by $(\Phi^{\OO_n},
\Phi^{\cN_n^{-1}\cM_n^{-1}})$, where $\cM_n = \fP_1\cdots\fP_g
\cap \OO_n$. We combine $\delta$ with $\pi$ to obtain the new
modular parametrization $\pi' : X_0(\fm\fn)\rightarrow E$ of
\S\ref{SectDegen}. Let $R \subset G_0$ be a set of representatives
for $G_0/G_1$. Notice that
\[
y_n = \Tr_{G_0}(\pi(x_n)) = \Tr_{R}\left(\sum_{\sigma\in
G_1}\pi(x_n)^{\sigma}\right) = \sum_{\sigma\in
R}\pi'(x'_n)^{\sigma}.
\]
For each $m\in\NN$ we choose
$\theta_m\in\Gal(K[\fp^{\infty}]/K[\fp^m])$ such that none of the
elements $\theta_m\sigma$ is geometric, for $\sigma\in R$. Then we
consider the composite map
\[
\begin{array}{ccccccc}
f_m : \NN & \stackrel{\rho}{\longrightarrow} &
X_0(\fm\fn)^{|R|}\times X_0(\fm\fn)^{|R|} &
\stackrel{\pi'}{\longrightarrow} &
E^{|R|}\times E^{|R|} \\
& & & \stackrel{\Sigma}{\longrightarrow} & E\times E &
\stackrel{\partial}{\longrightarrow} & E\quad\\
& & & & & & \\
 n & \longmapsto & ({x'_n}^{\sigma},
{x'_n}^{\theta_m\sigma})_{\sigma\in R} & \longmapsto &
(\pi'(x'_n)^{\sigma}, \pi'(x'_n)^{\theta_m\sigma})_{\sigma\in R}
\hspace{-0.8cm}
\\
& & & \longmapsto & (y_n, y_n^{\theta_m}) & \longmapsto & y_n-
y_n^{\theta_m}.
\end{array}
\]
We will show that $\rho$ is dominant, hence $f_m : \NN \rightarrow
E; \quad f_m(n) = y_n-y_n^{\theta_m}$ has finite fibres. In
particular, $f_m^{-1}(0)$ is finite, and the proof will be
complete.

Here the Andr\'e-Oort conjecture (see e.g. \cite{EdixPrep}) enters
the picture, for which we state the following characteristic-$p$
analogue (see \cite{BreuerPrep}).

\begin{conjecture}[Andr\'e-Oort]\label{AO}
Let $X=X_1\times\cdots\times X_n$ be a product of Drinfeld modular
curves, and let $Z\subset X$ be an irreducible algebraic
subvariety for which every projection $Z \rightarrow X_i$ is
dominant. Suppose $Z$ contains a Zariski-dense set of CM points.
Then $Z$ is a ``modular'' subvariety, which means the following.
There exist $g_1,\ldots,g_n\in \GL(k)$ and a partition
$\{1,\ldots,n\} = \coprod_{j=1}^m S_j$ such that $Z =
\prod_{j=1}^m Z_j$, and each $Z_j(\Cinf)\subset \prod_{i\in
S_j}X_i(\Cinf)$ is the image of $\Omega^*$ under the map $z\mapsto
([g_i(z)])_{i\in S_j}$.
\end{conjecture}

Let $I\subset\NN$ be an infinite subset. We suppose for the moment
that Conjecture \ref{AO} holds. In our case $X=X_0(\fm\fn)^{2|R|}$
and we take $Z$ to be a positive-dimensional irreducible component
of the Zariski-closure of $\rho(I)$ in $X$. Then Conjecture
\ref{AO} implies that either $Z=X$, or there exist some
$\sigma_i\neq\sigma_j \in R\cup\theta_m R$ such that the
projection of $Z$ onto the factor $X_0(\fm\fn)^2$ indexed by
$(\sigma_i,\sigma_j)$ is contained in some Hecke correspondence
$T_{\fd}$. But the latter case is impossible, as this would mean
that $\sigma_i\sigma_j^{-1}$ is geometric, contrary to the
definitions of $R$ and $\theta_m$. Thus $Z=X$, and we see that
$\rho$ is dominant. The result follows. More precisely, we have
shown:

\begin{theorem}\label{implication}
Suppose Conjecture \ref{AO} holds in the case where
$X=X_0(\fm\fn)^m$ and $Z\subset X$ contains a Zariski-dense set of
CM points $w=(w_1,\ldots,w_m)$ for which $\End(w_1) = \cdots =
\End(w_m) = \OO_n$ for some $n\in\NN$. Suppose further that the
map $\pi': X_0(\fm\fn)\rightarrow E$ is surjective. Then for every
infinite subset $I\subset\NN$ the group generated by
$\{y_n\;|\;n\in I\}$ in $E(H[\fp^{\infty}])$ has infinite rank.
\end{theorem}

In the special case where $k=\Fq(T)$, $A=\Fq[T]$ and $q$ is odd,
Conjecture \ref{AO} is known \cite{BreuerPrep}, and $\pi'$ is
surjective (Proposition \ref{nonconstant}), and so Theorem
\ref{MainTheorem} follows. \Endproof

{\bf Remark.} Theorems \ref{MainTheorem} and \ref{implication}
have two very easy mild generalizations: Firstly, fix an ideal
$\fc \subset A$ prime to $\fn$. Then we can construct Heegner
points corresponding to the orders $\OO_{\fc,n} = A +
\fc\fp^n\OO_K$, for which similar results hold. Secondly, let
$\chi : G_0 \rightarrow \{\pm 1\}$ be a character, and let
$K_{\chi}[\fp^{\infty}]$ be the subfield of $K[\fp^{\infty}]$
corresponding to $\ker(\chi)$. Note that $K_{\chi}[\fp^{\infty}] =
H[\fp^{\infty}]$ if $\chi$ is trivial. Then we may replace the
$G_0$-trace on $E$ by the $(G_0,\chi)$-trace
\[
\Tr_{G_0,\chi} : E(K[\fp^{\infty}]) \longrightarrow
E(K_{\chi}[\fp^{\infty}]);\quad x \longmapsto \sum_{\sigma\in G_0}
\chi(\sigma)x^{\sigma},
\]
and obtain a similar result. We leave the details to the dedicated
reader.


\end{document}